\documentclass[12pt]{article}
\usepackage{setspace}
\usepackage{amsmath,amssymb,amsthm}
\usepackage{mathrsfs}
\usepackage{geometry}
\usepackage{enumitem}
\usepackage{environ}
\usepackage{setspace}
\usepackage{titlesec}
\usepackage[all,cmtip]{xy}
\usepackage{tikz}
\usetikzlibrary{matrix,arrows}

\NewEnviron{prf}[1][]{\begin{proof}[\bf #1Proof]\BODY\end{proof}}{}
\NewEnviron{slt}[1][]{\begin{proof}[\bf #1	解]\BODY\end{proof}}{}
\newtheorem{definition}{Definition}[section]
\newtheorem{theorem}[definition]{Theorem}
\newtheorem{lemma}[definition]{Lemma}
\newtheorem{proposition}[definition]{Proposition}

\newtheorem{corollary}[definition]{Corollary}

\newcommand\MSC[1]{\textbf{MSC(2020)}: #1}
\newcommand\address[1]{#1}
\newcommand\email[1]{\emph{Email address}: #1}
  
\setlist[enumerate,1]{label=(\roman*)}
\setlist[enumerate,2]{label=(\alph*)}

\title{\textbf{Dynamical Mordell--Lang problem for automorphisms of surfaces in positive characteristic}}
\author{Junyi Xie\quad and\quad She Yang}
\date{}

\begin{document}
\begin{spacing}{1.25}

\maketitle

\begin{abstract}
We solve the dynamical Mordell--Lang problem in positive characteristic for automorphisms of projective surfaces.
\end{abstract}


\hspace{5.6em}height

\MSC{Primary: 37P55, 37P30; Secondary: 37P15.}

\section{Introduction}
In this paper, we work over an algebraically closed field $K$ of characteristic $p>0$. As a matter of convention, every variety is assumed to be integral but the closed subvarieties can be reducible. A surface is a two-dimensional variety. We denote $\mathbb{N}=\mathbb{Z}_{+}\cup\{0\}$. An arithmetic progression is a set of the form $\{mk+l|\ k\in\mathbb{Z}\}$ for some $m,l\in\mathbb{Z}$ and an arithmetic progression in $\mathbb{N}$ is a set of the form $\{mk+l|\ k\in\mathbb{N}\}$ for some $m,l\in\mathbb{N}$.

The dynamical Mordell--Lang conjecture is one of the core problems in the field of arithmetic dynamics. It asserts that for any endomorphism $f$ of a variety $X$ over $\mathbb{C}$, the return set $\{n\in\mathbb{N}|\ f^{n}(x)\in V(\mathbb{C})\}$ is a finite union of arithmetic progressions in $\mathbb{N}$ where $x\in X(\mathbb{C})$ is a point and $V\subseteq X$ is a closed subvariety. There is an extensive literature on various cases of this 0-DML conjecture (``0" stands for the characteristic of the base field). Two significant cases are as follows.
\begin{enumerate}
\item
If $X$ is a quasi-projective variety over $\mathbb{C}$ and $f$ is an \'etale endomorphism of $X$, then the 0-DML conjecture holds for $(X,f)$. See \cite{Bel06} and \cite[Theorem 1.3]{BGT10}.
\item
If $X=\mathbb{A}_{\mathbb{C}}^{2}$ and $f$ is an endomorphism of $X$, then the 0-DML conjecture holds for $(X,f)$. See \cite{Xie17} and \cite[Theorem 3.2]{Xie}.
\end{enumerate}

One can consult \cite{BGT16,Xie} and the references in there for more known results.

The statement of the 0-DML conjecture fails when the base field has positive characteristic. See \cite[Example 3.4.5.1]{BGT16} for an example. Indeed, the return set can be very complicated in positive characteristic. See \cite[Section 5]{XYb} and \cite[Section 5]{XY}. The $p$DML problem is known to be \emph{very} hard. It is proved in \cite{CGSZ21} that the $p$DML problem for endomorphisms of tori is \emph{equivalent} to solving some hard Diophantine equations.

For references toward the $p$DML problem, one can consult \cite{CGSZ21}, \cite[Theorem 1.4, Theorem 1.5]{Xie23}, \cite{Yang24}, and \cite{XYb,XY}.

In this article, we solve the $p$DML problem for automorphisms of surfaces. Our main theorem is as follows. Indeed, the statement also holds for automorphisms of proper surfaces and no change is needed in the proof. So we will just focus on the projective case.

\begin{theorem}
Let $X$ be a projective surface and let $f$ be an automorphism of $X$. Let $V\subseteq X$ be a closed subvariety and let $x\in X(K)$ be a point. Then the return set $\{n\in\mathbb{Z}|\ f^{n}(x)\in V(K)\}$ is a finite union of arithmetic progressions along with finitely many sets of the form $\{\frac{c_0+c_1q^m}{q-1}|\ m\in\mathbb{N}\}$, in which $q$ is a power of $p$ and $c_0,c_1$ are integers satisfying $q-1\mid c_0+c_1$. Moreover, the ``$p$-sets" will not appear unless $f$ is of bounded-degree.
\end{theorem}

We will completely determine the possible form of the return sets in Section 6. See Theorem 6.1.

In view of \cite[Corollary 1.3(ii)]{XYb}, we can prove the following corollary. We refer to \cite[Subsection 2.1]{XYb} for the definition of cohomologically hyperbolic endomorphisms and some relative knowledge.

\begin{corollary}
Let $X$ be a projective threefold and let $f$ be a cohomologically hyperbolic automorphism of $X$. Let $V\subseteq X$ be a closed subcurve and let $x\in X(K)$ be a point. Then the return set $\{n\in\mathbb{Z}|\ f^{n}(x)\in V(K)\}$ has the form described in Theorem 1.1.
\end{corollary}

Let $(X,f)$ be the dynamical system in Theorem 1.1. The proof of Theorem 1.1 is by dealing with three types of $(X,f)$ respectively. These three types are classified as follows.

Let $L\in\mathrm{Pic}(X)$ be a big and nef line bundle. We let $\lambda_1(f)=\lim\limits_{n\rightarrow\infty}((f^{n})^{*}L\cdot L)^{\frac{1}{n}}$. The limit does exist and this quantity is called the \emph{dynamical degree} of $f$, which is crucial in the study of arithmetic dynamics. See \cite{DS05,Dang20,Tru20} and \cite[Section 2.1]{Xie23} for the general theory of dynamical degrees. We say the automorphism $f$ is
\begin{enumerate}
\item
of bounded-degree (or elliptic), if the sequence $\{((f^{n})^{*}L\cdot L)|\ n\in\mathbb{N}\}$ is bounded;

\item
parabolic, if $\lambda_1(f)=1$ and the sequence $\{((f^{n})^{*}L\cdot L)|\ n\in\mathbb{N}\}$ is unbounded;

\item
hyperbolic, if $\lambda_1(f)>1$.
\end{enumerate}

We know that both the quantity $\lambda_1(f)$ and the type of $f$ are independent of the choice of $L$. Hence they are invariant under the semi-conjugation by generically finite surjective maps. See the references mentioned above.

Indeed, Theorem 1.1 is already proved in the literature if $f$ lies in the case (i) or case (iii). See \cite[Theorem 1.5,\ Remark 4.9]{XY}, \cite[Theorem 1.4]{Xie23}, and \cite[Corollary 1.3(ii)]{XYb}. Therefore, it turns out that we only need to prove Theorem 1.3 below.

\begin{theorem}
Let $X$ be a smooth projective surface and let $f$ be a parabolic automorphism of $X$. Let $V\subseteq X$ be a closed subvariety and let $x\in X(K)$ be a point. Then the return set $\{n\in\mathbb{Z}|\ f^{n}(x)\in V(K)\}$ is a finite union of arithmetic progressions.
\end{theorem}

The majority of this paper is dedicated to the proof of Theorem 1.3. We will use a height argument, which is of the same philosophy as that in \cite{XYb}. The key point of making the height argument work is to find \emph{two different} speeds of growth. But as we have mentioned in \cite[Section 5]{XYb}, in general it is hard to find two different speeds of growth in a dynamical system which is of zero entropy. To overcome this difficulty, we need to make use of Gizatullin's theorem for parabolic automorphisms of surfaces. This theorem asserts that the parabolic automorphisms will preserve a fibration. Then we can use this fibration to find two different speeds of growth. All of the height sequences in here will have a polynomial growth, but the degrees will be different.

The structure of this paper is as follows. We do some preparations in Section 2. We shall introduce the two main tools of this paper: the height machinery and the Gizatullin theorem for parabolic automorphisms. Then in Section 3, we deal with the case of abelian surfaces. In Section 4, we use the height machinery to cope with certain cases. Then we finish the proof of the main theorem in Section 5 by using an Albanese argument. Finally, we discuss the converse problem in Section 6. We will completely determine the form of return sets of surface automorphisms.

\section{Preparations}

In subsection 2.1, we will briefly recall Weil's height machinery. In subsection 2.2, we introduce Gizatullin's theorem for parabolic automorphisms of smooth projective surfaces, which says that such an automorphism preserves a fibration. Finally, in subsection 2.3, we state some easy facts about the dynamical Mordell--Lang problem for surfaces.

\subsection{The Weil height machinery}

For a detailed treatment of Weil's height machinery, we refer to \cite{Ser97}. For the purpose of our article, we will just list some basic facts here. We recommend the reader to take a look at \cite[Subsection 2.2]{XYb}. We inherit the setting in there and most of the contents in this subsection also appear in there.

In the following statements, we let $k$ be a field of characteristic $p>0$ equipped with a product formula. Then there is a standard logarithmic height function on $k$, or more generally, on the projective space $\mathbb{P}^{N}(k)$. We say that $k$ satisfies the \emph{Northcott property} if $\{x\in k|\ h(x)\leq A\}$ is a finite set for every $A>0$.

\begin{theorem}(\cite[Section 2.8]{Ser97})
Let $X$ be a projective variety over $k$. Denote $H$ as the quotient of the vector space of real-valued functions on $X(k)$ by the space of bounded functions on $X(k)$. Then there is a unique group homomorphism $L\mapsto h_{L}$ of $\mathrm{Pic}(X)$ to $H$ such that for every morphism $\phi:X\rightarrow\mathbb{P}_{k}^{N}$, we have $h_{\phi^{*}\mathcal{O}(1)}=h_{\phi}+O(1)$ in which $h_{\phi}(x)=h(\phi(x))$ is the naive height calculated on the projective space.
\end{theorem}

The following statements are immediate from definition.

\begin{lemma}
\begin{enumerate}
\item
Let $f:X\rightarrow Y$ be a morphism of projective varieties over $k$. Then $h_{f^*L}(x)=h_{L}(f(x))+O(1)$ as functions on $X(k)$ for every $L\in\mathrm{Pic}(Y)$.
\item
Let $X$ be a projective variety over $k$ and let $L$ be an ample line bundle on $X$. Then $h_L$ is bounded below. Suppose further that $k$ satisfies the Northcott property, then $\{x\in X(k)|\ h_{L}(x)\leq M\}$ is a finite set for every $M>0$ (and every representative of $h_{L}$).
\end{enumerate}
\end{lemma}

The following statement was proved in \cite[Proposition 2.9]{XYb}.

\begin{lemma}
Let $k$ be a finitely generated field extension of $\mathbb{F}_p$ of positive transcendence degree. Then we can make $k$ into a product formula field which satisfies the Northcott property.
\end{lemma}

We also need a bound on height of points on $\mathbb{P}^1$ and elliptic curves as below. Since part (i) of the following lemma is easy and part (ii) is well-known, we omit the proof. One can learn them from the knowledge of \cite{Ser97}.

\begin{lemma}
Let $k$ be a product formula field which satisfies the Northcott property.
\begin{enumerate}
\item
Let $a,b\in k^{\times}$. Suppose $a$ is not a root of unity. Then there exist $C_1,C_2>0$ such that $C_1n-C_2\leq h(ba^n)\leq C_1n+C_2$ for every nonnegative integer $n$.

\item
Let $E$ be an elliptic curve over $k$ and let $P,Q\in E(k)$. Suppose that $P$ is a non-torsion point. Let $L\in\mathrm{Pic}(E)$ be an ample line bundle and let $h_L$ be a representative of the height function associated to $L$. Then there exist $C_1,C_2,C_3>0$ such that $C_1n^2-C_3\leq h_L(nP+Q)\leq C_2n^2+C_3$ for every nonnegative integer $n$.
\end{enumerate}
\end{lemma}

\subsection{The Gizatullin theorem}

The Gizatullin theorem is a crucial result about parabolic automorphisms of smooth projective surfaces. The original reference is \cite{Giz80}, which deals with the rational surfaces in arbitrary characteristic (although there is an assumption that $\text{char}(K)\notin\{2,3\}$, please see also \cite[Remark 4.3]{CD12}). It turns out that the rational case is the most important one, as other cases can be treated by basic arguments.

In this subsection, we use the good survey \cite{Gri16} as a reference. Notice that this reference works over $\mathbb{C}$ while we focus on the positive characteristic case. So we will follow the reference and provide short proofs at some time, instead of directly citing the result.

We start with some general facts proved by pure linear algebra method in \cite[Section 3]{Gri16}. This result is characteristic-free. We denote $\mathrm{N}^1(X)$ as the numerical class group of a projective variety $X$, which is a finite free $\mathbb{Z}$-module. We will use ``$\equiv$" to denote numerical equivalence. The following theorem can be learned from \cite[Proposition 4.1]{Gri16} and the proof of \cite[Proposition 3.1]{Gri16}.

\begin{theorem}
Let $X$ be a smooth projective surface and let $f$ be a parabolic automorphism of $X$.
\begin{enumerate}
\item
There exists a unique primitive $f^*$-invariant nef class $\theta\in\mathrm{N}^1(X)$. 

\item
For every $f^*$-invariant class $\alpha\in\mathrm{N}^1(X)_{\mathbb{R}}$, we have $\theta\cdot\alpha=0$. In particular, we have $\theta^2=\theta\cdot K_X=0$.

\item
There exists $\theta'\in\mathrm{N}^1(X)$ and a positive integer $M$ such that $f^*\theta'-\theta'=M\theta$.
\end{enumerate}
\end{theorem}

Since Theorem 2.5(iii) is an important ingredient of this paper, we provide a bit more detail about its connection with \cite[Proposition 3.1]{Gri16}. By considering the intersection quadratic form on $\mathrm{N}^1(X)_{\mathbb{R}}$ and the automorphism $f^*$ of $\mathrm{N}^1(X)_{\mathbb{R}}$, we fit into the pure linear algebra setting of loc. cit.. Furthermore, we lie in case (3) of loc. cit. because $f$ is assumed to be parabolic. Now in fact, the unique invariant isotropic ray in loc. cit. is $\mathbb{R}_{+}\theta$, and hence we can find $\theta'\in\mathrm{N}^1(X)_{\mathbb{R}}$ such that $f^*\theta'-\theta'=\theta$ as that ray directs the first vector of that nontrivial Jordan block in loc. cit.. We can further make $\theta'$ to have rational coefficients as $\theta\in\mathrm{N}^1(X)$ and $f^*$ preserves the lattice $\mathrm{N}^1(X)$.

The following proposition is the key step in the proof of Gizatullin's theorem. We shall follow the arguments in \cite[Section 4]{Gri16} to prove it by the classification of surfaces. We refer to \cite{Bad01} for the basic knowledge of classification of surfaces in arbitrary characteristic.

\begin{proposition}
Let $X,f$ and $\theta$ be as in Theorem 2.5. Then there exists a globally generated line bundle $L$ on $X$ and a positive integer $N$, such that $L\equiv N\theta$ and $h^0(L)\geq2$.
\end{proposition}

\begin{prf}
We will prove that there exist $L$ and $N$ which satisfy the last two conditions. Then one can modify $L$ to make it globally generated by adapting the argument in the proof of \cite[Proposition 4.3]{Gri16}.

As we have mentioned before, the proof is by using the classification of surfaces. We let $\kappa(X)$ be the Kodaira dimension of $X$.
\begin{enumerate}
\item
$\kappa(X)=-1$.

In this case, either $X$ is rational or the image of the Albanese map $X\rightarrow\mathrm{Alb}(X)$ is a curve (here we arbitrarily fix a base point on $X$ for the Albanese map). In the former case, the assertion is guaranteed by \cite[Theorem 4.1]{CD12}. So we only need to deal with the latter case. Let $C$ be the image of the Albanese map, which is a projective curve. Let $\pi:X\rightarrow C$ be the induced surjective morphism. Then the automorphism $f$ induces a (not necessarily group) automorphism of $\mathrm{Alb}(X)$, and hence induces an automorphism $g$ of $C$ which satisfies $\pi\circ f=g\circ\pi$. Let $H$ be an ample line bundle on $C$. Then the class of $\pi^*H$ in $\mathrm{N}^1(X)$ is both $f^*$-invariant and nef. So the class of $\pi^*H$ in $\mathrm{N}^1(X)$ is a positive multiple of $\theta$. Therefore, we can let $L$ be a suitable multiple of $\pi^*H$ and then the conditions are fulfilled.

\item
$\kappa(X)=0$.

In this case, we can assume that $X$ is minimal without loss of generality. See \cite[Lemma 4.7]{Gri16}, which is characteristic-free. Then there are four cases to consider.
\begin{enumerate}
\item
$X$ is a K3 surface. In this case, the assertion is guaranteed by the Riemann--Roch formula.

\item
$X$ is an Enriques surface. In this case, the assertion follows from the fact that any nef (and numerically nontrivial) line bundle on $X$ has a positive Iitaka dimension. This should be well-known --- see for example \cite[Proposition 2.4]{KKM20}. It can also be proved by using \cite[Theorem 7.11]{Bad01}.

\item
$X$ is an abelian surface. Then by \cite{naf}, we know that there exists a positive integer $M$ such that the numerical class $M\theta$ contains an effective line bundle. But then the assertion follows from the fact that any nontrivial effective line bundle on $X$ has a positive Iitaka dimension.

\item
$\mathrm{Alb}(X)$ is an elliptic curve. In this case, the assertion can be proved by the same argument as in case (i).
\end{enumerate}

\item
$\kappa(X)=1$.

See \cite[Lemma 4.6]{Gri16} for this case. Indeed, let $m$ be a positive integer such that $h^0(mK_X)\geq 2$. Then one can just take $L=mK_X-F$ where $F$ is the fixed part of $mK_X$.

\item
$\kappa(X)=2$.

This case will not come into the picture since $\mathrm{Aut}(X)$ is a finite group in this case. See \cite{MDLM78}.
\end{enumerate}
\end{prf}

Now we can prove Gizatullin's theorem.

\begin{theorem}(Gizatullin)
Let $X$ be a smooth projective surface and let $f$ be a parabolic automorphism of $X$. Then there exist a smooth projective curve $C$, a surjective morphism $\pi:X\rightarrow C$, and an automorphism $g$ of $C$, such that $\pi_*\mathcal{O}_X=\mathcal{O}_C$ and $\pi\circ f=g\circ\pi$.
\end{theorem}

\begin{prf}
Let $\theta,L,N$ be as in Theorem 2.5 and Proposition 2.6. Then $L$ induces a morphism $X\rightarrow\mathbb{P}(H^0(X,L))$, whose image is a curve since $L^2=0$. We let $\pi:X\rightarrow C$ be the Stein factorization of this morphism, then there exists an ample line bundle $H$ on $C$ such that $L=\pi^*H$. We will prove that there exists an automorphism $g$ of $C$ such that $\pi\circ f=g\circ\pi$.

Firstly, we notice that the numerical class of the fibers of $\pi$ is a positive multiple of $\theta$. Since $\theta$ is $f^*$-invariant, $\theta^2=0$, and $\pi$ has connected fibers, we get a bijective map $g:C(K)\rightarrow C(K)$ such that $\pi\circ f=g\circ\pi$ on $X(K)$. We need to prove that $g$ is indeed an automorphism of $C$.

Let $C_0'\subseteq C\times C$ be the image of the morphism $X\rightarrow C\times C$ given by $x\mapsto(\pi(x),\pi\circ f(x))$. Let $\pi_1',\pi_2':C_0'\rightarrow C$ be the two induced projections. Then $\pi_1'$ and $\pi_2'$ are bijective on closed points. Hence $C_0'$ is a projective curve, and we have $g=\pi_2'\circ\pi_1'^{-1}:C(K)\rightarrow C(K)$. Let $C_0$ be the normalization of $C_0'$ and let $\pi_1,\pi_2:C_0\rightarrow C$ be the morphisms induced by $\pi_1'$ and $\pi_2'$. Then the two sentences about $\pi_1'$ and $\pi_2'$ above also hold for $\pi_1$ and $\pi_2$. We can factor through $\pi_1$ as $C_0\stackrel{\mathrm{Frob}_{q_1}}\longrightarrow C_0^{(q_1)}\stackrel{\sim}\rightarrow C$ and $\pi_2$ as $C_0\stackrel{\mathrm{Frob}_{q_2}}\longrightarrow C_0^{(q_2)}\stackrel{\sim}\rightarrow C$ in which $q_1,q_2$ are powers of $p$. Then our task becomes to prove $q_1=q_2$.

Assume by contradiction that $q_1\neq q_2$. Without loss of generality, we assume that $q_1<q_2$. Then $g$ is an endomorphism of $C$ of degree $\frac{q_2}{q_1}$. But then by looking at the numerical class of the fibers of $\pi$, we deduce $f^{*}\theta=\frac{q_2}{q_1}\theta$ and get a contradiction. So $q_1=q_2$ and thus $g$ is an automorphism of $C$.
\end{prf}

\subsection{Dynamical Mordell--Lang problem for surfaces}

In this subsection, we state some easy facts about the dynamical Mordell--Lang problem for surfaces. The key point in these facts is that the dynamical Mordell--Lang problem is trivial for curves.

\begin{definition}
Let $X$ be a variety and let $f$ be an automorphism of $X$. We say that $(X,f)$ satisfies the DML property if for every closed subvariety $V\subseteq X$ and every point $x\in X(K)$, the return set $\{n\in\mathbb{Z}|\ f^{n}(x)\in V(K)\}$ is a finite union of arithmetic progressions.
\end{definition}

\begin{lemma}
Let $X$ be a surface and let $f$ be an automorphism of $X$. For a point $x\in X(K)$, we denote $\mathcal{O}_f(x)=\{f^n(x)|\ n\in\mathbb{Z}\}$ as the two-sided orbit of $x$.
\begin{enumerate}
\item
$(X,f)$ satisfies the DML property if and only if $(X,f^n)$ satisfies the DML property for some positive integer $n$.

\item
$(X,f)$ satisfies the DML property if and only if for every $x\in X(K)$ whose two-sided orbit is dense in $X$ and every irreducible closed subcurve $C\subseteq X$, the intersection $\mathcal{O}_f(x)\cap C(K)$ is finite.

\item
Let $Y$ be a surface and let $g$ be an automorphism of $Y$. Suppose that there exists a dominant morphism $\pi:X\rightarrow Y$ such that $\pi\circ f=g\circ\pi$. Then $(X,f)$ satisfies the DML property if and only if $(Y,g)$ satisfies the DML property.
\end{enumerate}
\end{lemma}

We omit the proof since the proof is routine. Notice Lemma 2.9(ii) guarantees that in order to tackle the dynamical Mordell--Lang problem for surfaces, we only need to focus on the dense orbits.

\section{The case of abelian surfaces}

In this section, we prove Theorem 1.3 in the case that $X$ is an abelian surface. It is a philosophy of Ghioca and Scanlon that a dynamical system in positive characteristic is likely to \emph{not} have the DML property if it ``comes from algebraic groups". So in some sense, this is the most ``dangerous" case of Theorem 1.3. Indeed, it turns out that a special treatment is needed for this case. We will prove this case by using the Gizatullin theorem and an explicit height argument. We remark that our automorphisms of abelian varieties do not need to send 0 to 0. If so, then we will say ``group automorphism".

\begin{proposition}
Let $X$ be an abelian surface and let $f$ be a parabolic automorphism of $X$. Then $(X,f)$ satisfies the DML property.
\end{proposition}

\begin{prf}
Firstly, Theorem 2.7 says that there exist a smooth projective curve $C$, a surjective morphism $\pi:X\rightarrow C$, and an automorphism $g$ of $C$, such that $\pi_*\mathcal{O}_X=\mathcal{O}_C$ and $\pi\circ f=g\circ\pi$. Using the proposition on page 84 of \cite{Mum08}, one can prove that a surjective morphism from an abelian variety to $\mathbb{P}^1$ cannot have connected fibers. So the genus $g(C)$ cannot be 0. If $g(C)>1$, then $g$ is of finite order and hence every orbit of $f$ is not dense in $X$. So according to Lemma 2.9(ii), we have nothing to prove. Therefore, we may assume that $C$ is an elliptic curve without loss of generality and change the notation from $C$ to $E$.

Now we forget the setting of Theorem 2.7 and just remember that we have an elliptic curve $E$, a surjective morphism $\pi:X\rightarrow E$ satisfying $\pi_*\mathcal{O}_X=\mathcal{O}_E$, and an automorphism $g$ of $E$ such that $\pi\circ f=g\circ\pi$. We want to prove that $(X,f)$ satisfies the DML property. By substituting $f$ by an appropriate iteration, we may assume that $g$ is a translation of $E$ (notice that $f$ remains parabolic after the iteration). Then by compositing $\pi$ with a certain translation of $E$, we may further assume that $\pi$ is a group homomorphism. We denote $P=f(0)\in X(K)$, then $g$ is the translation of $E$ by $\pi(P)$. We write $f=\tau_{P}\circ f_0$ where $\tau_P$ is the translation map by $P$ on $X$ and $f_0$ is a group automorphism of $X$. Then we have $\pi\circ f_0=\pi$.

Let $E_0=\mathrm{ker}(\pi)$. Since the general fibers of $\pi$ are integral, we know that $E_0$ is an elliptic curve. Since $f$ is parabolic, it cannot be a translation of $X$. Thus $f_0\neq1$ and therefore we have $E_0=\mathrm{Im}(f_0-1)$. So $f_0(E_0)=E_0$ and hence $f_0$ induces a group automorphism of $E_0$, which we denote as $f_1$. Also, we can find a homomorphism $s:E\rightarrow X$ such that the homomorphism $p:E_0\times E\rightarrow X$ given by $(e_0,e)\mapsto e_0+s(e)$ is an isogeny. Then we denote $f_2$ as the homomorphism $E\rightarrow E_0$ given by $e\mapsto(f_0-1)\circ s(e)$. This is well-defined as $E_0=\mathrm{Im}(f_0-1)$.

Now we let $F_0$ be the group automorphism of $E_0\times E$ given by $(e_0,e)\mapsto(f_1(e_0)+f_2(e),e)$. Then one can verify that $p\circ F_0=f_0\circ p$. Pick a point $(a_0,a)\in(E_0\times E)(K)$ such that it maps to $P$ under $p$. Then the automorphism $F=\tau_{(a_0,a)}\circ F_0$ of $E_0\times E$ satisfies $p\circ F=f\circ p$. Notice that $F$ is also a parabolic automorphism, and we only need to prove that $F$ satisfies the DML property according to Lemma 2.9(iii).

For every positive integer $n$, we calculate that $F^n$ maps $(e_0,e)\in(E_0\times E)(K)$ to the point $(f_1^n(e_0)+\sum\limits_{i=0}^{n-1}f_1^i(a_0+f_2(e))+\sum\limits_{j=1}^{n-1}jf_1^{n-1-j}(f_2(a)),e+na)$. Notice that $f_1$ has finite order as it is a group automorphism of $E_0$. Let $d$ be the order of $f_1$. If $d\geq2$, then we have $\sum\limits_{i=0}^{d-1}f_1^{i}=0$. Hence $F^d$ is a translation of $E_0\times E$, which contradicts with the fact that $F$ is parabolic. So we conclude that $f_1=\mathrm{id}_{E_0}$. Thus we may shorten the formula as $F^n(e_0,e)=(e_0+n(a_0+f_2(e))+\frac{n(n-1)}{2}f_2(a),e+na)$, which indeed holds for every integer $n$. Now for the same reason that $F$ is parabolic, we can see that $f_2:E\rightarrow E_0$ is nonzero. Hence $f_2$ is an isogeny.

As our goal is to prove that $F$ satisfies the DML property, we may assume that $a$ is non-torsion since otherwise there will be no dense orbit of $F$. Then $f_2(a)$ is also non-torsion because $f_2$ is an isogeny. Now Lemma 3.2 below completes the proof.
\end{prf}

\begin{lemma}
Let $E_1$ and $E_2$ be elliptic curves. Let $x,y,z\in E_1(K)$ and let $a,b\in E_2(K)$. Suppose that both $z$ and $b$ are non-torsion. Then for any irreducible closed subcurve $C\subseteq E_1\times E_2$, the set $\{n\in\mathbb{Z}|\ (x+ny+\frac{n(n-1)}{2}z,a+nb)\in C(K)\}$ is finite.
\end{lemma}

\begin{prf}
We only need to prove that $\{n\in\mathbb{N}|\ (x+ny+\frac{n(n-1)}{2}z,a+nb)\in C(K)\}$ is a finite set. We may assume that $C$ dominates both $E_1$ and $E_2$ since otherwise the proof is easy. We find a finitely generated subfield $k\subseteq K$ such that the elliptic curves $E_1$ and $E_2$, the points $x,y,z,a,b$, and the closed subvariety $C\subseteq E_1\times E_2$ are defined over $k$. By abusing notation, we do not change the name of these data. So we have $x,y,z\in E_1(k)$ and $a,b\in E_2(k)$, and $z$ and $b$ are still non-torsion. Moreover, the induced morphisms $p_1:C\rightarrow E_1$ and $p_2:C\rightarrow E_2$ are finite. We want to prove that the set $\{n\in\mathbb{N}|\ (x+ny+\frac{n(n-1)}{2}z,a+nb)\in C(k)\}$ is finite.

Since there are non-torsion points on elliptic curves over $k$, we see that $k$ is not a finite field. So according to Lemma 2.3, we can make $k$ into a product formula field which satisfies the Northcott property. Then we can apply the height machinery. Let $L_1\in\mathrm{Pic}(E_1)$ and $L_2\in\mathrm{Pic}(E_2)$ be ample line bundles. Then both $p_1^*L_1$ and $p_2^*L_2$ are ample line bundles on $C$. We let $h_{L_1}$ and $h_{L_2}$ be representatives of the corresponding height functions. Denote $P_n=x+ny+\frac{n(n-1)}{2}z$ and $Q_n=a+nb$ for simplicity. Then the following holds.
\begin{enumerate}
\item
By Lemma 2.4(ii), there exists $C_1,C_2>0$ such that $h_{L_2}(Q_n)\leq C_1n^2+C_2$ for every nonnegative integer $n$. Similarly, one can show that there exists $C_3,C_4>0$ such that $C_3n^4-C_4n^3\leq h_{L_1}(P_n)$ for every positive integer $n$.

\item
Let $M$ be a positive integer such that $Mp_2^*L_2-p_1^*L_1$ is an ample line bundle on $C$. Then by Lemma 2.2, there exists $C_5>0$ such that $h_{L_1}(p_1(x))\leq Mh_{L_2}(p_2(x))+C_5$ for every point $x\in C(k)$.
\end{enumerate}

Combining (i) and (ii) above, we deduce that the set $\{n\in\mathbb{N}|\ (P_n,Q_n)\in C(k)\}$ is finite and hence finish the proof.
\end{prf}

\section{A height argument}

In this section, we will use a height argument to prove Theorem 1.3 when the Albanese map $X\rightarrow\mathrm{Alb}(X)$ is surjective.

\begin{proposition}
Let $X$ be a projective surface and let $f$ be an automorphism of $X$. Let $C$ be a smooth projective curve, $g$ be an automorphism of $C$, and $\pi:X\rightarrow C$ be a surjective morphism satisfying $\pi\circ f=g\circ\pi$. Suppose that there exist a line bundle $L$ on $X$ and an ample line bundle $H$ on $C$ such that $f^*L-L=\pi^*H$. Then $(X,f)$ satisfies the DML property.
\end{proposition}

\begin{prf}
If the genus $g(C)>1$, then the automorphism $g$ will have a finite order because the automorphism group of $C$ is finite \cite{MDLM78}. Thus, we may assume that $g(C)\leq1$ since otherwise $f$ will have no dense orbits.

\textbf{Case 1:} $g(C)=0$.

In this case, we have $C\cong\mathbb{P}^1$. After a suitable conjugation and an appropriate iteration, we may assume that $C=\mathbb{P}^1$ and $g$ has the form $[u,v]\mapsto[u,av]$ for some $a\in K^{\times}$. This is because we are working over a field of positive characteristic. Moreover, the condition about the existence of line bundles remains valid after the iteration as $(f^n)^*L-L=\pi^*(\sum\limits_{i=0}^{n-1}(g^i)^*H)$ for every positive integer $n$. Also, we may assume that $a\in K^{\times}$ is not a root of unity.

According to Lemma 2.9(ii), it suffices to prove that $\{n\in\mathbb{Z}|\ f^n(x)\in C_0(K)\}$ is a finite set for every irreducible closed subcurve $C_0\subseteq X$ and every point $x\in X(K)$ for which $\pi(x)$ is not $[1,0]$ or $[0,1]$. We only need to prove that $\{n\in\mathbb{N}|\ f^n(x)\in C_0(K)\}$ is finite, because then we can apply this result to $f^{-1}$ (notice that $f^*L-L=(f^{-1})^*(-f^*L)-(-f^*L)$). Also, we may assume that $C_0$ dominates $\mathbb{P}^1$.

Now we can find a finitely generated subfield $k\subseteq K$ such that all of the data $X,f,\pi$, and $C_0\hookrightarrow X$ are defined over $k$. We may also let $x$ be a $k$-point of $X$ and let $a\in k^{\times}$. Moreover, we can assume that there exist line bundles on the models of $X$ and $\mathbb{P}^1$ which pullback to $L$ and $H$, respectively (of course, the model of $\mathbb{P}^1$ is $\mathbb{P}^1$). By abusing notation, we will not change the name of all these data. We still have $\pi\circ f=g\circ\pi$ where $g$ is the automorphism of $\mathbb{P}^1$ given by $[u,v]\mapsto[u,av]$. The line bundle $H$ is still ample and the equation $f^*L-L=\pi^*H$ is still valid. Also, the induced map $p:C_0\rightarrow\mathbb{P}^1$ is finite. We need to prove that the set $\{n\in\mathbb{N}|\ f^n(x)\in C_0(k)\}$ is finite.

Since $a\in k^{\times}$ is not a root of unity, we see that $k$ is not a finite field. By Lemma 2.3, we can make $k$ into a product formula field which satisfies the Northcott property. We let $h_L$ and $h_H$ be representatives of the corresponding height functions on $X(k)$ and $\mathbb{P}^1(k)$, respectively. As $f^*L-L=\pi^*H$, we know that there exists $C_1>0$ such that $|h_L(f(y))-h_L(y)-h_H(\pi(y))|\leq C_1$ for every $y\in X(K)$. Since $\pi(x)$ is not $[1,0]$ or $[0,1]$ and $a$ is not a root of unity, Lemma 2.4(i) shows that there exists $C_2,C_3>0$ such that $C_2n-C_3\leq h_H(g^n(\pi(x)))=h_H(\pi(f^n(x)))\leq C_2n+C_3$ for every nonnegative integer $n$. Combining these bounds, we see that $C_2n-C_4\leq h_L(f^{n+1}(x))-h_L(f^n(x))\leq C_2n+C_4$ for every $n\geq0$ where $C_4=C_1+C_3$.

Taking sum, we get $C_5,C_6>0$ such that $C_5n^2-C_6n\leq h_L(f^n(x))\leq C_5n^2+C_6n$ for every positive integer $n$. We argue that this speed of growth forces the set $\{n\in\mathbb{N}|\ f^n(x)\in C_0(k)\}$ to be finite. Let $i$ be the inclusion map $C_0\hookrightarrow X$, then $p=\pi\circ i$. Since $p$ is finite, we see that there exists a positive integer $M$ such that $Mp^*H-i^*L$ is ample. So by Lemma 2.2, there exists $C_7>0$ such that $Mh_H(p(y))\geq h_L(i(y))-C_7$ for every $y\in C_0(k)$. In particular, we have $Mh_H(\pi(f^n(x)))\geq h_L(f^n(x))-C_7$ if $f^n(x)\in C_0(k)$. Hence the assertion follows from the bounds above.

~

\textbf{Case 2:} $g(C)=1$.

In this case, the curve $C$ is an elliptic curve and hence we change the notation from $C$ into $E$. Since the proof in this case is very similar to the case above, we will just sketch the proof.

Firstly, by doing an appropriate iteration, we can assume that $g$ is the translation of $E$ by a non-torsion point $P\in E(K)$. Then we prove that $\{n\in\mathbb{Z}|\ f^n(x)\in C_0(K)\}$ is a finite set for every irreducible closed subcurve $C_0\subseteq X$ and every point $x\in X(K)$. As above, we only need to prove that $\{n\in\mathbb{N}|\ f^n(x)\in C_0(K)\}$ is finite, and we can assume that $C_0$ dominates $E$.

We find a finitely generated subfield $k\subseteq K$ such that all of the data involved here are defined over $k$. We make $k$ into a product formula field which satisfies the Northcott property, and apply the height machinery. The remaining procedure of the proof is just the same as above. We use the bound in Lemma 2.4(ii) and the key point is that the growth of $\sum\limits_{i=1}^{n}i^2$ is faster than the growth of $n^2$.
\end{prf}

Now let us recall the setting of Theorem 2.7. According to Theorem 2.5(iii) and the construction of $\pi:X\rightarrow C$ in the proof of Theorem 2.7, we can see that there exists a line bundle $L$ on $X$ and an ample line bundle $H$ on $C$ such that $f^*L-L\equiv\pi^*H$. We will use this observation to prove the following proposition. We often omit the base point when talking about the Albanese map. This should cause no ambiguity --- one can just arbitrarily fix a base point. We refer to the Appendix of \cite{Moc12} for the basic knowledge about Albanese varieties.

\begin{proposition}
Let $X$ be a smooth projective surface and let $f$ be a parabolic automorphism of $X$. Suppose that the Albanese map $X\rightarrow\mathrm{Alb}(X)$ is surjective, then $(X,f)$ satisfies the DML property.
\end{proposition}

\begin{prf}
The dimension of $\mathrm{Alb}(X)$ can be $0,1$, or 2.

If $\mathrm{Alb}(X)=0$, then $\mathrm{Pic}^0(X)=0$. Here we denote $\text{Pic}^0(X)\subseteq\text{Pic}(X)$ as the subgroup consists of all algebraically trivial line bundles. As a result, there exist a line bundle $L$ on $X$ and an ample line bundle $H$ on $C$ such that $f^*L-L=\pi^*H$ in the setting of Theorem 2.7 (see the paragraph above). Then the conclusion follows from Proposition 4.1.

Suppose that $\mathrm{Alb}(X)$ is an abelian surface. We denote $p$ as the Albanese map $X\rightarrow\mathrm{Alb}(X)$, which is surjective and generically finite. There exists a (unique) automorphism $g$ of $\mathrm{Alb}(X)$ such that $p\circ f=g\circ p$. Then $g$ is also a parabolic automorphism according to \cite[Theorem 1(ii)]{Dang20}. Hence the conclusion follows from Proposition 3.1 and Lemma 2.9(iii).

Now we assume that $\mathrm{Alb}(X)=E$ is an elliptic curve. Let $p:X\rightarrow E$ be the Albanese map. Same as above, there exists an automorphism $h$ of $E$ such that $p\circ f=h\circ p$. So the numerical class of the fibers of $p$ is both nef and $f^*$-invariant. Hence this class is a positive multiple of $\theta$, where $\theta$ is the class in Theorem 2.5(i). So according to Theorem 2.5(iii), we can find a line bundle $L$ on $X$ and an ample line bundle $H_1$ on $E$ such that $f^*L-L-p^*H_1\in\mathrm{Pic}^0(X)$. But by the construction of the Albanese variety (see for example \cite[Proposition A.6]{Moc12}), we know that $p$ induces an isomorphism $p^*:\mathrm{Pic}^0(E)\rightarrow\mathrm{Pic}^0(X)$. So there exists a line bundle $H_0$ of degree 0 on $E$ such that $f^*L-L=p^*(H_1+H_0)$. Then we can just take $H=H_1+H_0$ and all the hypotheses in Proposition 4.1 are fulfilled. Thus we finish the proof.
\end{prf}

\section{Finish of the proof}

In this section, we will use an Albanese argument to conclude the proofs of our main theorems. We start with a lemma, which is a special case of Ueno's theorem \cite[Theorem 10.9]{Ueno75} if the base field is $\mathbb{C}$. By tracking through the proof, one can see that the statement of the following lemma holds in arbitrary characteristic. We give a proof for completeness.

\begin{lemma}
Let $A$ be an abelian variety and let $X\subseteq A$ be an irreducible closed subvariety of dimension not greater than 2. Suppose that $\mathrm{Stab}_A(X)=\{a\in A(K)|\ a+X=X\}=\{0\}$. Then the automorphism group $\mathrm{Aut}(X)$ is finite.
\end{lemma}

\begin{prf}
We assume that $\mathrm{dim}(X)=2$ since otherwise the proof is easy. We may assume that the algebraically closed base field $K$ is uncountable since both the condition and the conclusion are stable under base extension.

Let $X_0$ be the minimal resolution of $X$. Its existence is proved by Lipman --- see paragraph A on page 155 of \cite{Lip78}. Then $X_0$ is a smooth projective surface. We prove that $X_0$ is a minimal surface. Let $X_0'$ be a relatively minimal model of $X_0$. Since every rational map from a nonsingular variety to an abelian variety is a morphism \cite[Theorem I.3.2]{Mil}, we see that the induced birational map $X_0'\dashrightarrow X$ is a morphism. As $X_0$ is the minimal resolution, we conclude that $X_0\cong X_0'$ and hence the assertion holds. Therefore, by the classification of surfaces \cite{Bad01}, we have either $\kappa(X_0)\geq1$ or $X_0$ is an abelian surface because $X_0$ has a maximal Albanese dimension.

If $X_0$ is an abelian surface, then $X$ is a translation of a two-dimensional abelian subvariety of $A$. This contradicts the hypothesis that $\mathrm{Stab}_A(X)=\{0\}$.

If $\kappa(X_0)=2$, then the automorphism group $\mathrm{Aut}(X_0)$ is finite \cite{MDLM78}. Thus the result holds as there is an injection $\mathrm{Aut}(X)\hookrightarrow\mathrm{Aut}(X_0)$.

Now we focus on the case $\kappa(X_0)=1$. We will deduce a contradiction in this case. Recall that such a surface admits an elliptic or quasi-elliptic fibration \cite[Theorem 9.9]{Bad01}. However, we can prove that $X_0$ cannot admit a quasi-elliptic fibration as follows. Firstly, we know that almost every fiber of a quasi-elliptic fibration is an integral rational curve. Also, there is no non-constant map from a rational curve to an abelian variety. So all these fibers must be contracted by the birational morphism $X_0\rightarrow X$, which is absurd. Thus $X_0$ admits an elliptic fibration.

Arguing as above, we can see that almost all fibers of this elliptic fibration will map to some elliptic curves in $A$ which are contained in $X$. But since there are only countably many one-dimensional abelian subvariety in $A$ \cite{LOZ96}, we see that there exists an elliptic curve $E\subseteq A$ such that there are infinitely many fibers whose image lies in the family $\{a+E|\ a\in A(K)\}$ (recall that we have assumed $K$ to be uncountable). The set of images of those infinitely many fibers is also infinite, and all of these images are contained in $X$. Therefore, we get $E(K)\subseteq\mathrm{Stab}_A(X)$ and hence deduce a contradiction. So we finish the proof.
\end{prf}

Now we can prove our main theorems.

\proof[Proof of Theorem 1.3]
Let $a:X\rightarrow A$ be the Albanese map of $X$, in which we abbreviate $\mathrm{Alb}(X)$ as $A$. Then there is a unique automorphism $g$ of $A$ such that $a\circ f=g\circ a$. Let $Y=\mathrm{Im}(a)$. Then $Y$ is an irreducible closed subvariety of $A$ and we denote $\mathrm{Stab}_A(Y)$ as the reduced closed subgroup of $A$ satisfying $\mathrm{Stab}_A(Y)(K)=\{a\in A(K)|\ a+Y=Y\}$. The notation in here and that in Lemma 5.1 should cause no confusion. Let $p:A\rightarrow B$ be the quotient map, where we abbreviate $A/\mathrm{Stab}_A(Y)$ as $B$. Let $Z=p(Y)$. Then $Z$ is an irreducible closed subvariety of $B$ satisfying $\mathrm{Stab}_B(Z)=\{0\}$.

Since $g(Y)=Y$, one can verify that there exists an automorphism $h$ of $B$ such that $p\circ g=h\circ p$. Denote $\pi=p\circ a$. Then we have $\pi\circ f=h\circ\pi$ and $Z=\mathrm{Im}(\pi)$. So $h(Z)=Z$ and hence $h$ induces an automorphism $u$ of $Z$. We have $\pi_0\circ f=u\circ\pi_0$ in which $\pi_0:X\rightarrow Z$ is the surjective morphism induced by $\pi$. According to Lemma 5.1, the automorphism group $\mathrm{Aut}(Z)$ is finite. Therefore, the automorphism $f$ will have no dense orbit unless $Z$ is a point. Hence we may assume that $Z$ is a point by taking Lemma 2.9(ii) into account.

Since $Z$ is point, we have $Y=\mathrm{Stab}_A(Y)$ (notice that $0\in Y(K)$). But since the Albanese map cannot factor through a proper abelian subvariety, we deduce that $Y=A$, i.e. the Albanese map $a:X\rightarrow A$ is surjective. Now the theorem follows from Proposition 4.2.
\endproof

\proof[Proof of Theorem 1.1]
Let $\pi:X_0\rightarrow X$ be the minimal resolution of $X$ (see the Introduction of \cite{Lip78}). Then $X_0$ is a smooth projective surface and $f$ induces an automorphism $f_0$ of $X_0$ which satisfies $\pi\circ f_0=f\circ\pi$. So we may assume that $X$ is smooth without loss of generality, as the types of $f_0$ and $f$ are same. Then the proof of Theorem 1.1 is immediate by combining \cite[Theorem 1.5,\ Remark 4.9]{XY}, \cite[Theorem 1.4]{Xie23}, and Theorem 1.3.
\endproof

\proof[Proof of Corollary 1.2]
If the two-sided orbit $\mathcal{O}_f(x)=\{f^n(x)|\ n\in\mathbb{Z}\}$ is dense in $X$, then \cite[Corollary 1.3(ii)]{XYb} guarantees that $\{n\in\mathbb{Z}|\ f^n(x)\in V(K)\}$ is a finite set. If $\mathcal{O}_f(x)$ is not dense, then the assertion follows from Theorem 1.1.
\endproof

\section{The converse problem}

In this section, we discuss the converse problem of Theorem 1.1. We will determine which sets can be realized as a return set of a surface automorphism.

The converse of $p$DML problem was studied in \cite{LN25}. However, due to the complicated examples constructed in \cite[Section 5]{XYb} and \cite[Section 5]{XY}, it seems hopeless to study the general converse of the $p$DML problem. Indeed, to our knowledge, the converse of the Skolem--Mahler--Lech problem in positive characteristic \cite[Conjecture 3.6]{Der07} is still open. Notice that in positive characteristic, the Skolem--Mahler--Lech problem is equivalent to the $p$DML problem for translations of tori.

Thanks to the relatively simple form of the return sets in Theorem 1.1, we can solve the converse problem in this case. The result is as follows.

\begin{theorem}
Let $X$ be a projective surface and let $f$ be an automorphism of $X$. Let $V\subseteq X$ be a proper closed subvariety and let $x\in X(K)$ be a point. Let $\mathcal{O}_f(x)=\{f^n(x)|\ n\in\mathbb{Z}\}$ be the two-sided orbit of $x$.
\begin{enumerate}
\item
If $\mathcal{O}_f(x)$ is not dense in $X$, then the return set $\{n\in\mathbb{Z}|\ f^n(x)\in V(K)\}$ is a finite union of arithmetic progressions.
\item
If $\mathcal{O}_f(x)$ is dense in $X$, then we can write the return set $\{n\in\mathbb{Z}|\ f^n(x)\in V(K)\}$ as the union of $\bigcup\limits_{i=1}^{s}\{\frac{c_{i}+d_{i}q^m}{q-1}|\ m\in\mathbb{N}\}$ with a finite set, in which $q$ is a power of $p$ and $c_{1},\dots,c_{s},d_{1},\dots,d_{s}$ are integers that satisfy $q-1\mid c_{i}+d_{i}$ and $q\nmid d_{i}$ for every $1\leq i\leq s$.
\end{enumerate}

Moreover, every set of the form mentioned above can be realized as a return set of a surface automorphism.
\end{theorem}

In short, this result says that the return sets of \emph{automorphisms} should be ``complete". For example, the return set cannot be $\{p^m|\ m\in\mathbb{Z}_+\}$. The flavor of this result is toward the bounded-degree automorphisms, since otherwise there will be no ``$p$-sets" according to Theorem 1.1.

We start with a lemma. The proof can be done by a direct calculation.

\begin{lemma}
Let $K=\overline{\mathbb{F}_p(t)}$. Let $q$ be a power of $p$ and let $d$ be a positive integer that is not a multiple of $q$. Let $\alpha\in K$ be a generator of the cyclic group $\mathbb{F}_{q}^{\times}$. Let $x$ be a $d$-th root of $t^d+\alpha$ in $K$. Then $\{n\in\mathbb{Z}|\ x^n=t^n+\alpha\}=\{dq^m|\ m\in\mathbb{N}\}$.
\end{lemma}

Now we can prove Theorem 6.1.

\proof[Proof of Theorem 6.1]
We will focus on case (ii) since both the proof part and the construction part of case (i) are easy. We just make a remark: the key point for the construction part of case (i) is that the ordinary elliptic curves have torsion points of any (exact) order.

Now we consider the construction part of case (ii). We let the base field $K=\overline{\mathbb{F}_p(t)}$ and consider the automorphisms of $\mathbb{P}^2$. Let $d\in\mathbb{Z}_+$ be the least common multiple of $d_1,\dots,d_s$ and write $d=d_1d_1'=\cdots=d_sd_s'$. Then $d$ is still not a multiple of $q$. We have $\bigcup\limits_{i=1}^{s}\{\frac{c_{i}+d_{i}q^m}{q-1}|\ m\in\mathbb{N}\}=\bigcup\limits_{i=1}^{s}\{\frac{d_i'c_{i}+dq^m}{d_i'(q-1)}|\ m\in\mathbb{N}\}$. We write $a_i=d_i'c_i$ and $b_i=d_i'(q-1)$ for $1\leq i\leq s$.

We embed $\mathbb{G}_m^2$ into $\mathbb{P}^2$ by the map $(x,y)\mapsto[x,y,1]$. We write $[n]$ as the multiple-by-$n$ map of $\mathbb{G}_m^2$. We let $\alpha\in K$ be a generator of the cyclic group $\mathbb{F}_{q}^{\times}$ and let $x_0$ be a $d$-th root of $t^d+\alpha$ in $K$, as in Lemma 6.2. Let $f$ be the automorphism of $\mathbb{P}^2$ given by the formula $[x,y,z]\mapsto[x_0x,ty,z]$. Let $C\subseteq\mathbb{G}_m^2$ be the line $x=y+\alpha$. For $1\leq i\leq s$, we let $C_i$ be the closure of $[b_i]^{-1}(a_i\cdot(x_0,t)+C)\subseteq\mathbb{G}_m^2$ in $\mathbb{P}^2$. Then we have that $\{n\in\mathbb{Z}|\ f^n([1,1,1])\in C_i(K)\}=\{\frac{a_i+dq^m}{b_i}|\ m\in\mathbb{N}\}$ for every $1\leq i\leq s$. Hence the construction part is done.

Now we turn to the proof of part (ii). If $f$ is not of bounded-degree, then Theorem 1.1 says that the return set is finite. So we may focus on the bounded-degree case. In this case, the proof is by tracking through the arguments in \cite{XY}. We will be concise and just concentrate on the key point.

By certain reduction steps, we reduce to prove the same assertion for a translation of a two-dimensional isotrivial semi-abelian variety, say $G$. Our task then becomes to describe the intersection set of a dense cyclic group $\Gamma=\mathbb{Z}\cdot g\subseteq G(K)$ with an irreducible closed subcurve $C\subseteq G$. We denote $[n]$ as the multiple-by-$n$ map of $G$. By the calculations in \cite[Proposition 2.11]{XY}, we see that the $p$-sets will not appear unless a power of the Frobenius coincides with a power of $[p]$ on an infinite subgroup of $\Gamma$. But since such an infinite subgroup is also dense, we see that those two maps are equal. In particular, the map $[p]$ must be injective if the $p$-sets occur (i.e. if the return set is infinite). So without loss of generality, we may assume that $[p]$ is injective. Hence so does $[q]$ in which $q$ is an arbitrary power of $p$. 

We only need to prove the following assertion: if $\{\frac{c+dq^m}{q-1}|\ m\in\mathbb{Z}_+\}\cdot g\subseteq C(K)$, then $\frac{c+d}{q-1}\cdot g\in C(K)$ as well. Here $q$ be a power of $p$ and $c,d$ are integers satisfying $q-1\mid c+d$. This assertion guarantees that the return set must be complete, as we expected.

To prove this, we may assume $d\neq0$. Then $C=\overline{\{\frac{c+dq^m}{q-1}|\ m\in\mathbb{Z}_+\}\cdot g}$. Also, we have $[q](C)-c\cdot g=\overline{\{\frac{c+dq^{m+1}}{q-1}|\ m\in\mathbb{Z}_+\}\cdot g}$. Hence $[q](C)-c\cdot g=C$ as $[q]$ is finite. But since $[q](\frac{c+d}{q-1}\cdot g)=(c+\frac{c+dq}{q-1})\cdot g\in c\cdot g+C(K)=[q](C)(K)$, we conclude that $\frac{c+d}{q-1}\cdot g\in C(K)$ as $[q]$ is injective. Thus we finish the proof.
\endproof

\section*{Acknowledgements}
We are grateful to Guoquan Gao and Xiangqian Yang for many useful discussions. We thank the anonymous referee for carefully reading the manuscript.

This work is supported by the National Natural Science Foundation of China Grant No. 12271007.

\bibliographystyle{alpha}
\bibliography{reference}

\address{Beijing International Center for Mathematical Research, Peking University, Beijing 100871, China}

\email{xiejunyi@bicmr.pku.edu.cn}

~

\address{Beijing International Center for Mathematical Research, Peking University, Beijing 100871, China}

\email{ys-yx@pku.edu.cn}

\end{spacing}
\end{document}